\documentclass[12pt]{amsart}
\usepackage{amssymb}
\usepackage{hyperref}

\newcommand{\nn}[1]{(\ref{#1})}
\newcommand{\eqn}[1]{equation~\nn{#1}}

\newtheorem{theorem}[subsection]{Theorem}
\newtheorem{lemma}[subsection]{Lemma}
\newtheorem{proposition}[subsection]{Proposition}
\newtheorem{definition}[subsection]{Definition}

\newtheorem{example}[subsection]{Example}

\newcommand{\sfrac}[2]{{\textstyle \frac{#1}{#2}}}

\newsavebox{\ttoiso}
\sbox{\ttoiso}{\begin{picture}(25,12)(-5,0)
\put(0,3){$\simeq$}
\put(0,-3){$\longrightarrow$}
\end{picture}}

\newcommand{\x}{\times}
\newcommand{\br}{\mbox{$\mathbb R$}}

\newcommand{\bh}{\mathbb H}
\newcommand{\bz}{\mathbb Z}
\newcommand{\bp}{\mathbb P}

\newcommand{\up}{\Upsilon}

\newcommand{\cg}{{\mathcal G}}

\newcommand{\cq}{{\mathcal Q}}

\newcommand{\ca}{{\mathcal A}}

\newcommand{\al}{\alpha}
\newcommand{\be}{\beta}
\newcommand{\ga}{\gamma}

\renewcommand{\phi}{\varphi}

\newcommand{\si}{\sigma}

\newcommand{\ta}{\tau}

\newcommand{\Ga}{\Gamma}
\newcommand{\La}{\Lambda}
\newcommand{\Om}{\Omega}

\newcommand{\fg}{\mathfrak g}

\newcommand{\cv}{{\mathcal V}}
\newcommand{\cw}{{\mathcal W}}

\def\sideremark#1{\ifvmode\leavevmode\fi\vadjust{\vbox to0pt{\vss
 \hbox to 0pt{\hskip\hsize\hskip1em
 \vbox{\hsize3cm\tiny\raggedright\pretolerance10000
 \noindent #1\hfill}\hss}\vbox to8pt{\vfil}\vss}}}%
\begin{document}

\title{Generalized planar curves and quaternionic geometry}

\author{Jaroslav Hrdina and Jan Slov\'ak}


\address{Department of Algebra and Geometry, Faculty of Science, 
Masaryk University in Brno, Jan\'a\v ckovo n\'am. 2a, 662 95 Brno, 
Czech Republic}

\begin{abstract}
Motivated by the analogies between the projective and the almost quaternionic
geometries, we first study the generalized planar curves and mappings. We
follow, recover, and extend the classical approach, see e.g. 
\cite{Mi1, Mi2}. Then we exploit the impact of the general results in the
almost quaternionic geometry. In particular we show, that the natural class
of $\bh$--planar curves coincides with the class of all geodesics of the so
called Weyl connections and preserving this class turns out to be the
necessary and sufficient condition on diffeomorphisms to become 
morphisms of almost quaternionic geometries.
\end{abstract}

\maketitle

\frenchspacing

Various concepts generalizing geodesics of affine connections have been
studied for almost quaternionic and similar geometries.  Let us point out
the generalized geodesics defined via generalizations of normal coordinates,
cf. \cite{BE} and \cite{CS}, or more recent \cite{CSZ, Z}.  Another class of
curves was studied in \cite{Mi2} for the hypercomplex structures with
additional linear connections. The latter authors called a curve $c$
quaternionic planar if the parallel transport of each of its tangent vectors
$\dot c(t_0)$ along $c$ was quaternionic colinear with the tangent field
$\dot c$ to the curve. Yet another natural class of curves is given by the
set of all unparameterized geodesics of the so called Weyl connections, i.e.
the connections compatible with the almost quaternionic structure with
normalized minimal torsion. The latter connections have remarkably similar
properties for all parabolic geometries, cf. \cite{CS}, and so their name
has been borrowed from the conformal case.  In the setting of almost
quaternionic structures, they were studied first in \cite{O} and so they are
also called Oproiu connections, see \cite{AMP}.

The first author showed in \cite{H} that actually the concept of
quaternionic planar curves was well defined for the almost quaternionic
geometries and their Weyl connections. Moreover, it did not depend on the
choice of a particular Weyl connection and it turned out that the
quaternionic planar curves were just all unparameterized geodesics of all
Weyl connections. 

The aim of this paper is to find further analogies of Mike\v s's classical
results in the realm of the almost quaternionic geometry. On the way we
simplify, recover, and extend results on generalized planar mappings,
explain results from \cite{H}, and
finally we show that morphisms of almost quaternionic geometries are just
those diffeomorphisms which leave invariant the class of all
unparameterized geodesics of Weyl connections. 

\noindent{\bf Acknowledgments.} 
The authors have been supported by GACR, grants Nr.
201/05/H005 and 201/05/2117, respectively.
The second author also gratefully acknowledges support from the Royal
Society of New Zealand via Marsden Grant no. 02-UOA-108 during 
writing essential part of this paper. The authors are also grateful to
Josef Mike\v s for numerous discussions and to Dmitri Alekseevsky for
helpful comments.

\section{Motivation and background on quaternionic geometry}\label{1}

There are many equivalent definitions of almost quaternionic geometry to be
found in the literature. Let us start with the following one:

\begin{definition}\label{1.1}
Let $M$ be a smooth manifold of dimension $4n$. 
An almost hypercomplex structure on $M$ is a
triple $(I,J,K)$ of smooth affinors in $\Ga(T^*M\otimes TM)$ satisfying
$$
I^2=J^2=-E,\quad K= I\circ J=-J\circ I
$$
where $E=\operatorname{id}_{TM}$. 

An almost quaternionic structure is a rank four subbundle
$\cq\subset T^*M\otimes TM$ 
locally generated by the identity $E$ and a hypercomplex
structure.
\end{definition}

An almost complex geometry on a $2m$--dimensional
manifold $M$ is given by the choice of the affinor $J$ satisfying $J^2=-E$.
Let us observe, that such a $J$ is uniquely determined within the rank two
subbundle $\langle E, J\rangle\subset TM$, up to its sign. Indeed, if
$\hat J=aE + bJ$, then the condition $\hat J^2=-E$ implies $a=0$ and
$b=\pm1$.

Thus we may view
the almost quaternionic geometry as a straightforward generalization of this
case. Here, a similar simple computation reveals that the rank three 
subbundle $\langle I, J, K \rangle$ is invariant of the choice of the
generators and this is the definition we may find in \cite{AMP}. More
explicitly, different choices will always satisfy $\hat I= aI+bJ+cK$ with
$a^2+b^2+d^2=1$, and similarly for $J$ and $K$.
Let us also remark that the 4--dimensional almost
quaternionic geometry coincides with 4-dimensional conformal Riemannian
geometries.

\subsection{The frame bundles}\label{1.2}
Equivalently, we can define an almost quaternionic structure $\cq$ on $M$ as
a reduction of the linear frame bundle $P^1M$ to an appropriate structure group,
i.e. as a G--structure with the structure group of all automorphisms
preserving the subbundle $\cq$. We may view such frames as linear 
mappings $T_xM\to \bh^n$ 
which carry over the multiplications by $i,j,k\in\bh$ onto some of the
possible choices for $I,J,K$.  Thus, a 
further reduction to a fixed hypercomplex
structure leads to the structure group $GL(n,\bh)$
of all quaternionic linear mappings on $\bh^n$.
Additionally, we have to allow morphisms which do not leave the 
affinors $I$, $J$, $K$ invariant but change them within the subbundle $\cq$.
As well known, the resulting group is 
\begin{equation}\label{G0}
G_0= GL(n,\bh)\x_{\bz_2}Sp(1)
\end{equation}
where $Sp(1)$ are the unit quaternions in $GL(1,\bh)$, see e.g. \cite{S}. 

We shall write $\cg_0\subset P^1M$ for this principal $G_0$--bundle defining
our structure.

The simplest example of such a structure is well understood as the
homogeneous space 
$$
\bp_n\bh = G/P
$$
where $G_0\subset P$ are the subgroups in $G=SL(n+1,\bh)$
\begin{align}
G_0 &= \left\{
\begin{pmatrix}
a & 0
\\
0 & A
\end{pmatrix};\ A\in GL(n,\bh),\ \operatorname{Re}(a\operatorname{det}A)=1 
\right\},
\\
P& = \left\{
\begin{pmatrix}
a & Z
\\
0 & A
\end{pmatrix};\  
\begin{pmatrix}
a & 0
\\
0 & A
\end{pmatrix} \in G_0, 
\ Z \in (\bh^n)^*
\right\}
.\end{align}

Since $P$ is a parabolic subgroup in the semisimple Lie group $G$, 
the almost quaternionic geometry is an instance of the so called parabolic
geometries. All these geometries enjoy a rich and quite uniform theory
similar to the classical development of the conformal Riemannian and
projective geometries, but
we shall not need much of this here. We refer the reader to \cite{CS} and
the references therein.

\subsection{Weyl connections} \label{1.3}
The classical prolongation procedure for G--structures starts with finding a
minimal available torsion for a connection belonging to the structure on 
the given manifold $M$. Unlike the projective and conformal Riemannian
structures where torsion free connections always exist, the torsion has to
be allowed for the almost quaternionic structures in general in dimensions
bigger than four. The standard
normalization comes from the general theory of parabolic geometries and we
shall not need this in the sequel. The details may be found for example in
\cite{GS}, \cite{BE}, for another and more classical point of view 
see \cite{S}. The only essential point for us is that
all connections compatible with the given geometry sharing the unique
normalized torsion are parameterized by smooth one--forms on $M$. In analogy
to the conformal Riemannian geometry we call them {\em Weyl connections} for
the given almost quaternionic geometry on $M$.

The almost quaternionic geometries with Weyl connections without torsion 
are called {\em quaternionic geometries}.

From the point of view of prolongations of G--structures, 
the class of all Weyl connections defines a reduction $\cg$ of the semiholonomic
second order frame bundle over the manifold $M$ to a principal subbundle 
with the structure group $P$, while the individual Weyl connections
represent further reductions of $\cg$ to principal subbundles $\cg_0\subset
\cg$ with the structure group $G_0\subset P$.  
Thus, the Weyl connections $\nabla$ are in bijective correspondence with
$G_0$--equivariant sections $\si:\cg_0\to\cg$ of the natural projection.   

As mentioned above, the difference of two Weyl connections is a one--form and, 
also in full analogy to the conformal geometry, there are neat formulae for the
change of the covariant derivatives of two such connections $\hat\nabla$ and
$\nabla$ in terms of their difference $\up=\hat\nabla - \nabla \in
\Om^1(M)$. 

\subsection{Adjoint tractors}\label{1.4}
In order to understand the latter formulae, we introduce the so called
adjoint tractors. They are sections of the vector bundle (called usually
{\em adjoint tractor bundle}) over $M$
$$
\ca = \cg_0\x_{G_0} \fg
,$$
see (1) through (3) for the definition of the groups.
The Lie algebra $\fg={\mathfrak sl}(n+1,\bh)$ carries the $G_0$--invariant 
grading 
\begin{equation}\label{grading}
\fg=\fg_{-1}\oplus \fg_0\oplus \fg_1
\end{equation}
where 
\begin{align*}
\fg_0 &= \left\{
\begin{pmatrix}
a & 0
\\
0 & A
\end{pmatrix};\ A\in {\mathfrak gl}(n,\bh),\ a\in \bh
\ \operatorname{Re}(a+\operatorname{Tr}A)=0 
\right\},
\\
\fg_1& = \left\{
\begin{pmatrix}
0 & Z
\\
0 & 0
\end{pmatrix};\   
Z \in (\bh^n)^*
\right\},\quad
\fg_{-1} = \left\{
\begin{pmatrix}
0 & 0
\\
X & 0
\end{pmatrix};\  
X \in \bh^n
\right\}
.\end{align*}
Moreover, $TM = \cg_0\x_{G_0} \fg_{-1}$, $T^*M=\cg_0\x_{G_0}\fg_1$, 
and we obtain on the level of vector bundles
$$
\ca = \ca_{-1}\oplus \ca_0\oplus \ca_1 = TM \oplus \ca_0\oplus T^*M
$$
where $\ca_0= \cg_0\x \fg_0$ is the adjoint bundle of the Lie algebra
$\fg_0$.
 
The key feature of $\ca$ is that 
all further $G_0$--invariant objects on $\fg$ are carried over to the
adjoint tractors, too. In particular, the Lie bracket on $G$ induces an
algebraic bracket $\{\ ,\ \}$ on $\ca$ such that each fibre has the
structure of the graded Lie algebra (\ref{grading}).

Now we may write down easily the transformation formula. 
Let $\hat\nabla$ and $\nabla$ be two Weyl connections,
$\hat\nabla-\nabla=\up\in \Ga(\ca_1)$. Then for all tangent vector fields $X,Y\in
\Ga(\ca_{-1})$, i.e. sections of the $\fg_{-1}$-component of the adjoint
tractor bundle, 
\begin{equation}\label{e0}
\hat\nabla_X Y = \nabla _X Y + \{\{X,\up\},Y\} = \nabla_XY + \up(X)Y,
\end{equation}
where $\up(X)\in \fg_0 = \frak s\frak p(1) +  \frak g\frak l(n, \bh)$, 
see \cite{CS} or \cite {BE, GS} for the proof. 
Notice that the internal bracket results in an endomorphism on $TM$, while
the external bracket is just the evaluation of this endomorphism on $Y$ (all
this is read off the brackets in the Lie algebra easily).  

\section{Generalized planar curves and mappings}\label{2}

Various geometric structures on manifolds are defined as smooth distributions 
in the vector bundle $T^*M\otimes TM$ of all endomorphisms of the tangent
bundle. We have seen the two examples of almost complex and almost
quaternionic structures above. Let us extract some formal
properties from these examples.

\begin{definition}\label{2.1}
Let $A$ be a smooth $\ell$--dimensional vector subbundle in $T^*M\otimes TM$,
such that the identity affinor $E= \operatorname{id}_{TM}$ restricted to
$T_xM$ belongs to $A_xM\subset T^*_xM\otimes T_xM$ at each point $x\in M$.
We say that $M$ is equipped by an $A$--structure.  
\end{definition}    

For any tangent vector $X\in T_xM$ we shall write $A(X)$ for the vector 
subspace 
$$
A(X) = \{F(X);\ F\in A_xM\}\subset T_xM
$$
and we call $A(X)$ the {\em $A$--hull of the vector $X$}. Similarly,
the {\em $A$--hull of a vector field} will be the subbundle in $TM$
obtained pointwise. Notice that the dimension of such a subbundle in $TM$ may
vary pointwise. For every smooth parameterized curve $c:\br \to M$ we write
$\dot c$ and $A(\dot c)$ for the tangent vector field and its $A$-hull along
the curve $c$. 

For any vector space $V$, we say that a vector subspace 
$A\subset V^*\otimes V$ of automorphisms
is of {\em generic rank $\ell$} 
if the dimension of $A$ is $\ell$, and
the subset of vectors $(X, Y)\in V\x V$, such that the $A$--hulls 
$A(X)$ and $A(Y)$ generate a vector subspace $A(X)\oplus A(Y)$ 
of dimension $2\ell$, is open and dense. An $A$--structure is said to be of
generic rank $\ell$ if $A_xM$ has this property for each point $x\in M$.
 
Let us point out some examples: 
\begin{itemize}
\item The $\langle E \rangle$--structure is of
generic rank one on all manifolds of dimensions at least 2. 
\item Any almost
complex structure or almost product structure $\langle E, J \rangle$ is of
generic rank two on all manifolds of dimensions at least 4. 
\item Any almost
quaternionic structure is of generic rank four on all manifolds of
dimensions at least 8.
\end{itemize}
Similar conclusions apply to paracomplex and parahermition structures.

\begin{definition}\label{2.2}
Let $M$ be a smooth manifold with a given $A$--structure and a linear
connection $\nabla$. A smooth curve $c:\br\to M$ is told to be $A$--planar if 
$$
\nabla_{\dot c}{\dot c} \in A(\dot c)
.$$   
\end{definition}

Clearly, $A$ planarity means that the parallel transport of any tangent vector
to $c$ has to stay within the $A$--hull $A(\dot c)$ of the 
tangent vector field $\dot c$ along the curve. Moreover, this concept does
not depend on the parametrization of the curve $c$.

\begin{definition}\label{2.3}
Let $M$ be a manifold with a linear connection $\nabla$ and an
$A$--structure, while $N$ be
another manifold with a linear connection $\hat\nabla$ and a $B$--structure. 
A diffeomorphism $f:M\to N$ is called
$(A,B)$--planar if each $A$--planar curve $c$ on $M$ is mapped onto the
$B$--planar curve $F\circ c$ on $N$.  
\end{definition}

\begin{example}\label{2.4}
The 1--dimensional $A=\langle E\rangle$ structure must be given just as the 
linear hull of the identity affinor $E$, by the
definition. Obviously, the $\langle E\rangle$--planar curves on a manifold $M$ with a linear
connection $\nabla$ are exactly the unparameterized geodesics. Moreover,
two connections $\nabla$ and $\bar\nabla$ 
without torsion are projectively equivalent (i.e. they share the same
unparameterized geodesics) if and only if their
difference satisfies $\bar\nabla_X Y - \nabla_X Y = \al(X)Y + \al(Y)X$ for some
one--form $\al$ on $M$. The latter condition can be rewritten as
\begin{equation}
\bar\nabla - \nabla \in \Ga(T^*M \odot \langle E\rangle)\subset \Ga(S^2T^*M\otimes TM) 
\end{equation}
where the symbol $\odot$ stays for the symmetrized tensor product.
\end{example}

The latter condition on projective structures may be also rephrased in the
terms of morphisms: A diffeomorphism $f:M\to M$ is called geodesical 
(or an automorphism of the projective structure) if
$f\circ c$ is an unparameterized geodesic for each geodesic $c$ and this
happens if and only if the symmetrization of the difference $f^*\nabla -
\nabla$ is a section of $T^*M\odot \langle E\rangle$. 
We are going to generalize the above example in the rest of this section. 

In the case $A=\langle E \rangle$, the $(\langle E \rangle, B)$--planar
mappings are called simply {\em $B$--planar}. They map each geodesic curve 
on $(M, \nabla)$ onto a $B$--planar curve on $(N,\hat \nabla, B)$. 

Each $\ell$ dimensional $A$ structure $A\subset T^*M\otimes TM$ determines
the distribution $A^{(1)}$ in $S^{2}T^*M\otimes TM$, given at any point
$x\in M$ by
$$
A^{(1)}_xM = \{ \al_1\odot F_1+\dots+ \al_\ell\odot F_\ell;\ \al_i\in T^*_xM,
F_i\in A_xM\}
.$$

\begin{theorem}\label{2.5}
Let $M$ be a manifold with a linear connection $\nabla$, let $N$ be a
manifold of the same dimension 
with a linear connection $\hat \nabla$ and an $A$--structure of
generic rank $\ell$, and suppose $\operatorname{dim}M\ge 2\ell$. Then a
diffeomorphism $f:M\to N$ is $A$--planar if and only if
\begin{equation}\label{e2}
\operatorname{Sym}(f^*\hat\nabla-\nabla)\in f^*(A^{(1)})
\end{equation}
where $\operatorname{Sym}$ denotes the symmetrization of the difference of
the two connections.
\end{theorem}

The proof is based on a purely algebraical lemma below. Let us first
observe that the entire claim of the theorem is of local character.
Thus, identifying the objects on $N$ with their pullbacks on $M$, we may
assume that $M=N$ and $f=\operatorname{id}_M$. 

Next, let us observe that the $A$--planarity of $f:M\to N$ does not at all
depend on the possible torsions of the connection. Indeed, we always test
expressions of the type $\nabla_{\dot c}\dot c$ for a curve $c$ and thus
deforming $\nabla$ into $\nabla'=\nabla+T$ by adding some torsion will not
effect the results. Thus, without any loss of generality, we may assume that
the connections $\nabla$ and $\hat\nabla$ share the same torsion, and then
we may omit the symmetrization from \eqn{e2}.    

Finally, we may fix some (local) basis $E=F_0, F_i$, $i=1,\dots\ell-1$, of $A$, 
i.e. $A=\langle F_0,\dots,F_{\ell-1} \rangle$. Then the condition in 
\eqn{e2} says
\begin{equation}\label{e3}
\hat\nabla = \nabla + \sum_{i=0}^{\ell-1} \al_i\odot F_i
\end{equation}
for some suitable one--forms $\al_i$ on $M$. Of course, the existence of
such forms does not depend on our choice of the basis of $A$.

The quite simplified statement we now have to prove is:

\begin{proposition}\label{2.5a} 
Let $M$ be a manifold of dimension at least $2\ell$, 
$\nabla$ and $\hat\nabla$ two connections on $M$ with
the same torsion, and consider an $A$--structure of generic rank $\ell$ on $M$. 
Then each geodesic curve with respect to $\nabla$ is $A$--planar with
respect to $\hat\nabla$ if and only if there are one--forms $\al_i$
satisfying \eqn{e3}. 
\end{proposition}

Assume first we have such forms $\al_i$, and let $c$ be a geodesic for
$\nabla$. 
Then \eqn{e3} implies $\hat\nabla_{\dot c}\dot c \in A(\dot c)$ so that $c$
is $A$--planar, by definition.

The other implication is the more difficult one. Assume each (unparameterized)
geodesic $c$ is $A$--planar. This implies that the symmetric difference tensor
$P=\hat\nabla-\nabla\in \Ga(S^2T^*M\otimes TM)$ satisfies
$$
P(\dot c,\dot c)= \hat\nabla_{\dot c}\dot c \in \langle\dot c, F_1(\dot
c),\dots, F_{\ell-1}(\dot c) \rangle
.$$ 
In fact, the main
argument of the entire proof boils down to a purely algebraic claim:

\begin{lemma}\label{2.5b} Let $A\subset V^*\otimes V$ be a vector subspace of generic
rank $\ell$, and assume that $P(X,X)\in A(X)$ for 
some fixed symmetric tensor $P\in V^*\otimes V^*\otimes V$ and 
each vector $X\in V$. 
Then the induced mapping $P:V\to V^*\otimes V$ has values in $A$. 
\end{lemma}

\begin{proof}
Let us fix a basis $F_0, F_1,\dots, F_{\ell-1}$ of $A$. Since $A$ is of
generic rank $\ell$, there is 
the open and dense subspace $\cv\subset V$ of all vectors $X\in V$
for which $\{F_0(X), F_1(X), \dots, F_{\ell-1}(X)\}$ are 
linearly independent. 
Now, for each $X\in \cv$ there are the unique coefficients 
$\al_i(X)\in\br$ such that
\begin{equation}\label{e4}
P(X,X) = \sum_{i=0}^{\ell -1}\al_i(X)F_i(X)
.\end{equation}
The essential technical step in the proof of our Lemma is to show that all
functions $\al_i$ are in fact restrictions of one--forms on $V$. Let
us notice, that $P$ is a symmetric bilinear mapping and thus it is 
determined by the restriction of $P(X,X)$ to arbitrarily small open
non--empty subset of the arguments $X$ in $V$. 

\noindent
{\bf Claim 1.} {\em If a symmetric tensor $P\in V^*\otimes V^*\otimes V$
is determined over the
above defined subspace $\cv$ by \nn{e4}, 
then the functions $\al_i:\cv\to \br$ are smooth and their
restrictions to the individual rays (half--lines) generated by vectors in $\cv$
are linear.}

Let us fix a local smooth basis $e_i\in V$, the dual basis $e^i$, and
consider the induced dual bases $e_I$ and $e^I$ on the multivectors and
exterior forms. Let us consider the smooth mapping 
$$
\chi:\La^\ell V\setminus \{0\}\to \La^\ell V^*,\quad 
\chi(\mbox{$\sum$} a_Ie^I)= \sum \frac{a_I}{\sum a_I^2}e_I
.$$ 
Now,
for all non--zero tensors 
$\Xi = \sum a_Ie^I$, the evaluation 
$\langle \Xi, \chi(\Xi)\rangle$ is the constant function 1, 
while $\chi(k\cdot \Xi)=k^{-1}\chi(\Xi)$. 

Next, we define for each $X\in\cv$
$$
\ta(X)= \chi\bigl(X\wedge F_1(X)\wedge\dots\wedge F_{\ell-1}(X)\bigr)
$$
and we may compute the unique coefficients $\al_i$ from \nn{e4}:
\begin{align*}
\al_0(X) &= \langle P(X,X)\wedge F_1(X)\wedge F_2(X)\wedge\dots\wedge F_{\ell-1}(X), \ta(X) 
\rangle 
\\
\al_1(X) &= \langle X\wedge P(X,X)\wedge F_2(X)\wedge\dots\wedge F_{\ell-1}(X), \ta(X) 
\rangle 
\\
&\vdots
\\
\al_{\ell-1}(X) &= \langle X\wedge F_1(X)\wedge F_2(X)
\wedge\dots\wedge P(X,X), \ta(X)  \rangle 
.\end{align*}
In particular, this proves the first part our Claim 1.

Let us now consider a fixed vector $X\in \cv$. The defining formula \nn{e4}
for $\al_i$ implies  $\al_i(kX)=k\al_i(X)$, for each real number $k\ne0$.
Passing to zero with positive $k$ shows that $\al$ does have the limit $0$
in the origin and so we may extend the definition of $\al_i$'s (and validity
of formula \nn{e4}) to the entire
cone $\cv\cup \{0\}$ by setting $\al_i(0)=0$ for all $i$. 

Finally, along the ray 
$\{tX; t>0\}\subset \cv$, the derivative 
$\sfrac d{dt}\al(tX)$ has the constant value $\al(X)$.
This proves the rest of Claim 1.   

Now, in order to complete the proof of Lemma \ref{2.5b}, we have to prove
the following assertion.

\noindent
{\bf Claim 2.} {\em If a symmetric tensor $P$ is determined over the
above defined subspace $\cv\cup \{0\}$ by \nn{e4}, then 
the coefficients $\al_i$ are linear one--forms on $V$ and
the tensor $P$ is given by 
$$
P(X,Y) = \frac12 \sum_{i=0}^{\ell-1}\bigl( \al_i(Y)F_i(X) 
+ \al_i(X)F_i(Y)\bigr)
.$$
}

The entire tensor $P$ is obtained through polarization from its evaluations
$P(X,X)$, $X\in TM$, 
\begin{equation}\label{e5}
P(X,Y)= \sfrac12\bigl( P(X+Y,X+Y) - P(X,X) - P(Y,Y) \bigr),
\end{equation}
and again, the entire tensor is determined by its values on arbitrarily
small non--empty open subset of $X$ and $Y$ in each fiber. 

The summands on the right hand side have values in the following subspaces:
\begin{align*}
P(X+Y,X+Y) \in\ &\langle X+Y, F_1(X+Y),\dots, F_{\ell-1}(X+Y) \rangle\subset
\\
&\langle X, F_1(X),\dots
F_{\ell-1}(X), Y, F_1(Y),\dots, F_{\ell-1}(Y) \rangle,   
\\
P(X,X) \in\ &\langle X, F_1(X),\dots F_{\ell-1}(X) \rangle
,\\
P(Y,Y) \in\ &\langle Y, F_1(Y),\dots F_{\ell-1}(Y) \rangle.
\end{align*}
Since we have assumed that $A$ has generic rank $\ell$, the
subspace $\cw\in V\x_M V$
of vectors $(X,Y)$ such that all the values 
$$
\{ X, F_1(X),\dots,
F_{\ell-1}(X), Y, F_1(Y),\dots, F_{\ell-1}(Y)\}
$$ 
are linearly independent is
open and dense. Clearly $\cw\subset \cv\x\cv$. Moreover, if
$(X,Y)\in \cw$ than $F_0(X+Y),\dots, F_{\ell-1}(X+Y)$ are independent,
i.e. $X+Y\in \cv$. Inserting \nn{e4} into
\nn{e5}, we obtain
$$
P(X,Y)=\sum_{i=0}^{\ell-1}\bigl(d_i(X,Y)F_i(X) + 
e_i(X,Y) F_i(Y) \bigr)
.$$
For all $(X,Y)\in \cw$, 
the coefficients $d_i(X,Y)=\sfrac12(\al_i(X+Y)-\al_i(X))$ at $F_i(X)$, and 
$e_i(X,Y)= \sfrac12(\al_i(X+Y)-\al_i(Y))$ at $F_i(Y)$ in the
latter expression are uniquely determined. The symmetry of $P$ implies
$d_i(X,Y)= e_i(Y,X)$. If $(X,Y)\in\cw$ then also $(sX,tY)\in \cw$ for all
non-zero reals $s,t$ and the linearity of $P$ in the individual arguments
yields for all real parameters $s,t$
$$
std_i(X,Y) = sd_i(sX,tY)
.$$ 
Thus the functions $\al_i$ satisfy
$$
\al_i(sX+tY)-\al_i(sX) = t(\al_i(X+Y)-\al_i(X))
.$$
Since $\al_i(tX)=t\al_i(X)$, in the limit $s\to 0$ this means 
$$
\al_i(Y) = \al_i(X+Y)-\al_i(X) 
.$$
Thus $\al_i$ are additive over the open and dense set $(X, Y)\in\cw$.
Choosing a basis of $V$ such that each couple of basis elements is in
$\cw$, this shows that $\al_i$ are restrictions of linear forms, as
required. 
\end{proof}

Now the completion of the proof of Theorem \ref{2.5} is straightforward.
Folowing the equivalent local claim in the Proposition \ref{2.5a} and the
pointwise algebraic description of $P$ achieved in Lemma \ref{2.5b}, we just
have to apply the latter Lemma to individual fibers over the points $x\in M$
and verify, that the linear forms $\al_i$ may be chosen in a smooth way. But
this is obvious from the explicit expression for the coefficients $\al_i$ in
the proof of Claim 1 above.  

\begin{theorem}\label{2.6}
Let $M$ be a manifold with a linear connection $\nabla$ and an $A$--structure, 
$N$ be a manifold of the same dimension 
with a linear connection $\hat \nabla$ 
and a $B$--structure with generic rank $\ell$. Then a
diffeomorphism $f:M\to N$ is $(A,B)$--planar if and only if $f$ is
$B$--planar and $A(X)\subset (f^*(B))(X)$ for all $X\in TM$.
\end{theorem}
\begin{proof}
As in the proof of Theorem \ref{2.5}, we may restrict ourselves to some open
submanifolds, fix generators $F_i$ for $B$, assume that 
$f=\operatorname{id}_M$
and both connections $\nabla$ and $\hat\nabla$ share the same torsion,
and prove the equivalent local assertion to our theorem:

\noindent
{\bf Claim.} {\em Each $A$--planar curve $c$ with respect to $\nabla$ 
is $B$--planar with respect to $\hat\nabla$, if and only if the
symmetric difference tensor $P=\hat\nabla-\nabla$ is of the form \nn{e4}
with smooth one--forms $\al_i$, $i=0,\dots,\ell-1$ and $A(X)\subset B(X)$
for each $X\in TM$.}

Obviously, the condition in this statement is sufficient. So let us deal
with its necessity. 

Since every $(A,B)$--planar mapping is also 
$B$--planar, Theorem \ref{2.5} (or the equivalent Proposition in its proof) 
says that 
$$
P(X,X)=\sum_{j=0}^\ell \al_i(X)F_i(X)
$$
for uniquely given smooth one--forms $\al_i$. 

Now, consider a fixed $F\in A$ and suppose $F(X)\notin B(X)$. 
Since we assume that all $\langle E, F
\rangle$--planar curves $c$ in $M$ are
$B$--planar, we may proceed exactly as 
in the beginning of the proof of Theorem \ref{2.5}
to deduce that 
$$
P(X,X)=\sum_{j=0}^\ell \al_i(X)F_i(X) + \be(X)F(X)
$$
on a neighborhood of $X$, with some unique functions $\al_i$ and $\be$.

The comparison of the latter two unique expressions for $P(X,X)$ shows that
$\be(X)$ vanishes. But since $F(X)\ne X$, there definitely are curves which
are $\langle E,F \rangle$--planar and tangent to $X$, 
but not $\langle E \rangle$--planar. 
Thus, the assumption in the theorem would lead to $\be(X)\ne0$. Consequently, 
our choice $F(X)\notin B(X)$ cannot be achieved and
we have proved $A(X)\subset B(X)$ for all $X\in TM$.
\end{proof}

\section{Results on quaternionic geometries}\label{3}

The main result of this section is:

\begin{theorem}\label{3.1}
Let $f:M\to M'$ be a diffeomorphism between two almost quaternionic
manifolds of dimension at least eight. 
Then $f$ is a morphism of the geometries if and only if it
preserves the class of unparameterized geodesics of all Weyl connections on
$M$ and $M'$.
\end{theorem}

This theorem will follow easily from 
the results of Section \ref{2} and its proof 
requires only a few  quite simple formal steps. On the way we shall also
give a complete description of all geodesics of the Weyl connections in
terms of the $Q$--planar curves where $Q$ is the almost quaternionic
structure.

\begin{lemma}\label{3.2}
A curve $c:\br\to M$ is $Q$--planar with respect to at least
one Weyl connection $\nabla$ on $M$ if and only if $c$ is 
$Q$--planar with respect to
all Weyl connections on $M$.
\end{lemma}

\begin{proof}
For a Weyl connection
$\nabla$ and a curve $c:\br\to M$, the defining equation for $Q$--planarity
reads $\nabla_{\dot c}\dot c \in Q(\dot c)$. If we choose 
some hypercomplex structure within $Q$, we may rephrase this condition as:
$\nabla_{\dot c}\dot c = \dot c \cdot q$ where $\dot c(t)$ is a curve in the
tangent bundle $TM$ while $q(t)$ is a suitable 
curve in quaternions $\bh$. Now the
formula \nn{e0} for the deformation of the Weyl connections implies
$$
\hat\nabla_{\dot c}\dot c = \nabla_{\dot c}\dot c + \{\{\dot c,\up\},\dot c\}
= \nabla_{\dot c}\dot c + 2\dot c\cdot \up(\dot c).
$$
Indeed, this is the consequence of the computation of the Lie bracket in
$\fg$ of the corresponding elements $\dot c\in \fg_{-1}$, $\up\in \fg_1$:
\begin{align*}
[[\dot c,\up],\dot c] &\simeq \left[\left[
\begin{pmatrix} 0&0\\ \dot c&0\end{pmatrix},
\begin{pmatrix} 0&\up\\ \dot 0&0\end{pmatrix}
\right], \begin{pmatrix} 0&0\\ \dot c&0\end{pmatrix}
\right]
\\
&=
\begin{pmatrix} 0&0\\ 2\dot c \cdot \up(\dot c)&0\end{pmatrix}
\simeq 2\dot c \cdot \up(\dot c)
,\end{align*}
where $\up(\dot c)$ is the standard evaluation of the linear form $\up\in
\fg_1=(\bh^n)^*$ on the vector $\dot c\in \fg_{-1}=\bh^n$. 
Thus we see that if there is such a quaternion $q$ for one Weyl connection,
then it exists also for all of them. 
\end{proof}

\begin{definition}\label{3.3}
A curve $c:\br\to M$ is called $\bh$--planar if it is $Q$--planar with
respect to each Weyl connection $\nabla$ on $M$.
\end{definition}

\begin{theorem}\label{3.4}
Let $M$ be a manifold with an almost quaternionic geometry. Then a curve
$c:\br\to M$ is $\bh$--planar if and only if $c$ a geodesic of some Weyl
connection, up to parametrization.
\end{theorem}

\begin{proof}
Let us remark that $c$ is a geodesic for $\nabla$ if and only if $\nabla
_{\dot c}\dot c=0$. 
Thus, the statement follows immediately from the computation in the proof 
of Lemma \ref{3.2}. Indeed, if $c$ is $\bh$--planar, then choose any Weyl
connection $\nabla$ and pick up $\up$ so that $\hat\nabla_{\dot c}\dot c$
vanishes. 
\end{proof}

\subsection{Proof of Theorem \ref{3.1}}\label{3.5}
Every morphism of almost quaternionic geometries preserves the class of Weyl
connections and thus also the class of their geodesics.

We have to prove the oposit implication. This means, we have two manifolds
with almost quaternionic structures $(M,Q)$, $(N, Q')$ and a
diffeomorphism $f:M\to N$ which is $(Q,Q')$--planar. Then Theorem \ref{2.6}
implies that $Tf$ maps quaternionic lines in $T_xM$ to quaternionic lines in
$T_{f(x)}N$, i.e. $Q(X)=(f^*Q')(X)$ for each $X\in TM$ 
(since they both have the same dimension). In order to conclude the Theorem,
we have to verify $Q=f^*Q'$ instead, since this is exactly the requirement
that $f$ preserves the defining subbundles $Q$ and $Q'$.

Let us look at the 
subsets of all second jets of $\Bbb H$--planar curves. At each point, the
accelerations fill just the complete $Q$--hulls of the velocities (see
\cite{Z} for technicalities) and so,
for a given point $x\in M$ we may locally 
choose a smoothly parameterized system
$c_X$ of $\Bbb H$--planar curves with parameter $X\in T_xM$ such that 
$\dot c_X = X$ and
$\nabla_{\dot c_X}\dot c_X= \be(X) F(X)$ where $F$ is one of the generators
of $Q$ and $\be$ is a 1--form. Then 
$$
(\nabla_{\dot c_X}-\hat\nabla_{\dot c_X})\dot c_X = \be(X)F(X) +
\sum_k \al^k(X)F_k(X)
$$
where $F_k$ are the generators of $Q'$ and $\al^k$ are smooth 1--forms, cf.
the proof of Theorem \ref{2.5}. But this shows that, except of the zero set
of $\be$, $F(X) = \sum_k \ga^k(X)F_k(X)$ for some smooth $\ga^k$.
Since $F$ is linear in $X$, $\ga^k$ have to be constants and we are done.
\qed

\subsection{Final remarks}\label{3.6}
All curves in the four--dimensional quaternionic
geometries are $\bh$--planar by the definition. Thus this concept starts to
be interesting in higher dimensions only, and all of them are covered by
Theorem \ref{3.1}.

The class of the unparameterized geodesics of Weyl connections 
is well defined for all parabolic geometries.
Our result for the quaternionic geometries suggests the question, 
whether a similar statement holds for other geometries as well.

\end{document}